\documentstyle[12pt,leqno]{article}

\input boxedeps.tex
\SetepsfEPSFSpecial
\HideDisplacementBoxes

\def\reals{\hbox{\sl I\kern-.18em R \kern-.3em}}
\def\ints{\hbox{\sl Z\kern-.4em Z \kern-.3em}}
\def\nats{\hbox{\sl I\kern-.16em N \kern.05em}}
\def\rats{\hbox{\sl Q \kern-.83em\vrule height.59em depth0em \kern.87em}}
\def\complexes{\hbox{\sl\kern.50em I\kern-.50em C \kern.05em}}

\newtheorem{theorem}{Theorem} 
\newtheorem{proposition}{Proposition}

\def\np{\noindent}
\def\pf {\np {\bf Proof:} \ }
\def\endpf{$\|$ \bigskip}

\def\cH{{\cal H}}

\def\cK{{\cal K}}
\def\cL{{\cal L}}

\def\cX{{\cal X}}
\def\cY{{\cal Y}}
\def\cZ{{\cal Z}}

\def\bB{{\bf B}}
\def\bP{{\bf P}}

\def\bX{{\bf X}}
\def\bY{{\bf Y}}
\def\bZ{{\bf Z}}

\begin{document}
\title {Holonomic and Legendrian parametrizations of knots}
\author {Joan S. Birman \thanks{\np The first author acknowledges partial
support from the U.S.National Science Foundation under Grants DMS-9705019 
and DMS-9973232. The
second author is a graduate student in the Mathematics Department of Columbia
University. She was partially supported under the same grant, and also under
DMS-98-10750.}  
\\ jb@math.columbia.edu \and Nancy C. Wrinkle 
\\ wrinkle@math.columbia.edu}
\date{J.Knot Theory and its Ramifications {\bf 9} No. 3 (2000), p. 293-309.}
\maketitle
\centerline{Received 8 April 1999}
\begin{abstract} \noindent Holonomic parametrizations of knots 
were introduced in 1997 by Vassiliev, who proved that every knot type can be
given a holonomic parametrization. Our main result is that any two holonomic
knots which represent the same knot type are isotopic in the space of holonomic
knots. A second result emerges through the techniques used to prove the main
result: strong and unexpected connections between the topology of knots
and the algebraic solution to the conjugacy problem in the braid groups,
via the work of Garside.   We also discuss related parametrizations of Legendrian
knots, and uncover connections between the concepts of holonomic and
Legendrian parametrizations of knots.
\end{abstract}

\np \section{Introduction:} \label{section:Introduction} 
\begin{figure}[htpb]
\centerline{\BoxedEPSF{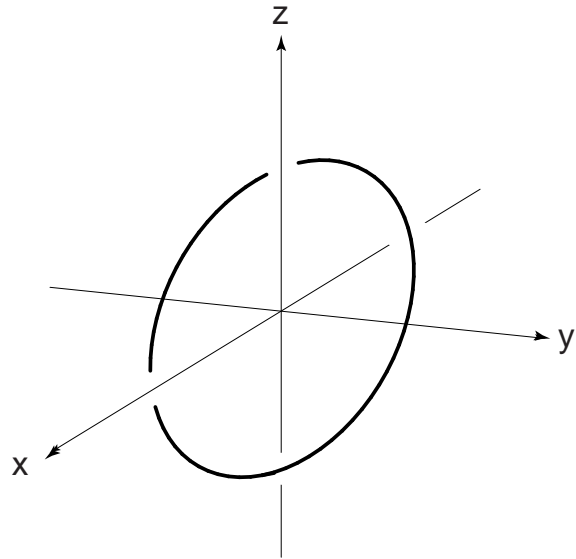 scaled 700}}
\caption{Three-space view of the unknot which is defined by $f(t) = cos(t).$}
\label{figure:holonomic unknot}
\end{figure}

\begin{figure}[htpb]
\centerline{\BoxedEPSF{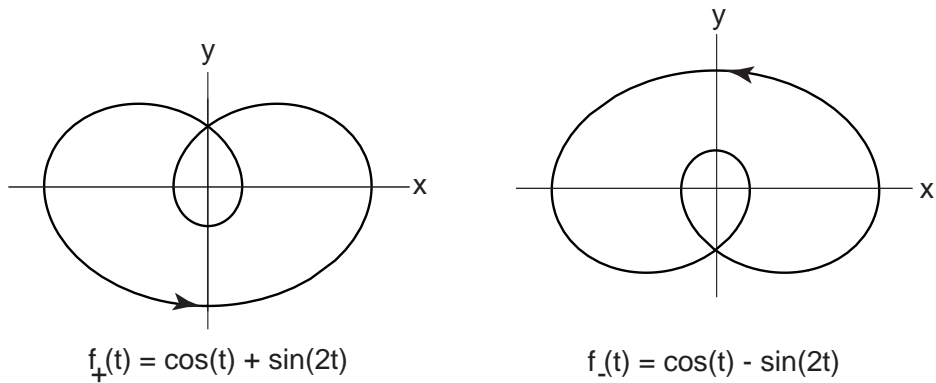 scaled 900}}
\caption{Two additional representations of the unknot}
\label{figure:2-braid unknots}
\end{figure}

Let $f:\reals \to
\reals$ be a
$C^\infty$ periodic function with period $2\pi$.   Following Vassiliev
\cite{Vass}, use
$f$ to define a map 
$\tilde{f}:S^1 \to \reals^3$ by setting $\tilde{f}(t) =
(-f(t), f'(t),-f''(t))$. Let $\pi$ be the restriction of $\tilde{f}$ to the
first two coordinates.  We call $\pi$ the {\em projection} of $K =
\tilde{f}(S^1)$  (onto the $xy$ plane). It turns out that with some
restrictions on the choice of the defining function $K$ will be a knot
and $\pi$ will yield a knot diagram with some very pleasant
properties, which we now begin to describe.  We highlight our assumptions about
$f$ with bullets, and describe their consequences:
\newpage
\begin{enumerate}
\item [(1)] We assume that $f$ is chosen so that $\tilde{f}(S^1)$ is 
a knot $K\subset\reals^3$, i.e. 
\begin{itemize}
\item There do not exist distinct points $t_1,t_2\in [0,2\pi)$ such
that \\
$(-f(t_1),f'(t_1),-f^{\prime\prime}(t_1)) =
(-f(t_2),f'(t_2),-f^{\prime\prime}(t_2)).$
\end{itemize}
Note that this implies that double points in the projection which are off the
$x$ axis are transverse. For, at a double point $(-f(t_1),f'(t_1))
= (-f(t_2),f'(t_2))$. The double point is transverse if the tangent vectors to
the projected image are distinct at $t_1$ and $t_2$. The tangent vectors are
$(-f'(t_1),f^{\prime\prime}(t_1))$ and $(-f'(t_2),f^{\prime\prime}(t_2))$. They
are distinct because $\tilde{f}$ is an embedding, which  implies
that  $f^{\prime\prime}(t_1) \neq f^{\prime\prime}(t_2)$.

\item [(2)]  The
reasoning used in (1) above shows that the tangent to $\pi(K)$ at an
instant when $\pi(K)$ crosses the $x$ axis is vertical. We observe that this
implies that if a double point occurred at an axis crossing, then it would
necessarily be a point where the two branches of
$\pi(K)$ had a common tangent. We rule out this behavior, which is not allowed
in a `regular' knot diagram,  by requiring that
$f$ be chosen so that all double points are away from the $x$ axis, i.e.
\begin{itemize}
\item If  $(-f(t_1),f'(t_1)) = (-f(t_2),f'(t_2))$ then $f'(t_1)\neq 0$.
\end{itemize}
\item [(3)] The singularities in a regular knot diagram are required to be at
most a finite number of transverse double points. To achieve that we need one
more assumption, i.e. 
\begin{itemize}
\item There do not
exist distinct points
$t_1,t_2,t_3\in [0,2\pi )$ such that  
$(-f(t_1),f'(t_1)) = (-f(t_2),f'(t_2)) = (-f(t_3),f'(t_3)) $.
\end{itemize} 
\end{enumerate}
Vassiliev observed that these conditions hold for generic $f$. A simple example is obtained by taking
$f(t) = cos(t)$,  giving the unknot which is pictured in Figure
\ref{figure:holonomic unknot}. 
Two additional examples are given in Figure \ref{figure:2-braid unknots},
which show the projections of the  knots which are
defined by the functions $f_\pm(t) = cos(t) \pm sin(2t)$.\\

\begin{figure}[htpb]
\centerline{\BoxedEPSF{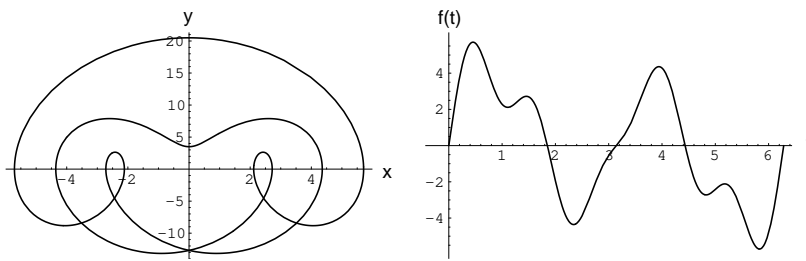 scaled 600}}
\caption{The function $f(t) = Sin(t) + 4Sin(2t) + Sin(4t) + 1.5Sin(5t)$ defines a holonomic trefoil
that winds continually anticlockwise but does not have a single point about which it winds.}
\label{figure:winding}
\end{figure}

\np There is an immediate suggestion of a closed braid in this parametrization,
for the following reasons. In the half-space $y > 0$, we know that $f'(t) > 0$,
so $-f(t)$ is decreasing.  Similarly, in the half-space $y < 0$, we have that
$-f(t)$ is increasing.   
Since $f$ is assumed to be generic, if
$\cK$ is crossing from the half-space $y < 0$ to the half-space $y > 0$ at $t_0$
then $z = -f''(t_0) < 0$ and if it is crossing from the half-space $y > 0$ to the
half-space $y < 0$ at $t_0$ then $z= -f''(t_0) > 0$.  Thus the projected image of
$K$ on the $xy$ plane winds  continually in an anticlockwise sense
(anticlockwise because the $x$ coordinate is $-f(t)$). The only reason it may
not already be a closed braid is that there may not be a single point on the $x$
axis which separates all of the axis-crossings with $f''(t) > 0$ from the
crossings where $f''(t) < 0$. An example of a holonomic knot which is not in
braid form is pictured in Figure \ref{figure:winding}.

\begin{enumerate}
\item[(4)] Analyzing the difficulty, we see that the number of zeros in
the $x$ coordinate, i.e. in the graph of
$f(t)$, is smaller than the number of zeros in the $y$ coordinate, i.e.
in the graph of $f'(t)$. We add one more requirement:
\begin{itemize}
\item The number of zeros in one cycle of $f$ is the same as the number of
zeros in one cycle of $f'$.
\end{itemize}
\end{enumerate}
When $f(t)$ is chosen to satisfy (1)-(4) our parametrization gives a closed
braid.  The braid index is then one-half the number of zeros in one cycle of $f$
(or of $f'$). When all these
conditions are satisfied our knot is said to have a  {\em holonomic
parametrization}.   
\\

\begin{figure}[htpb]
\centerline{\BoxedEPSF{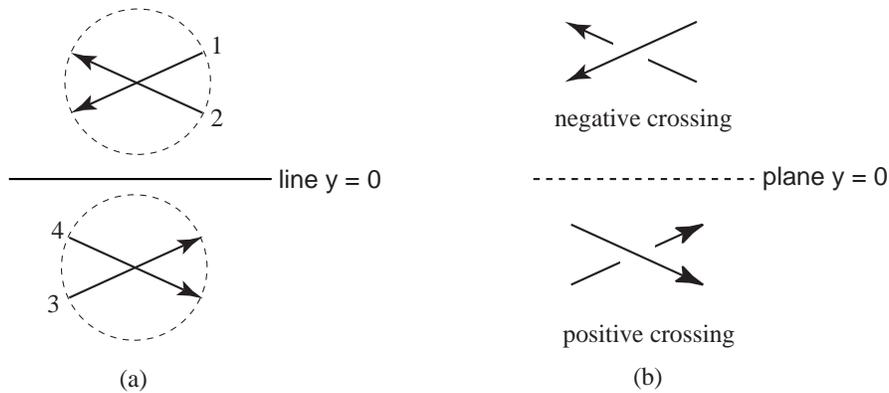 scaled 800}}
\caption{Determining the sign of a crossing: Figure (a) shows the
projected images onto the $xy$ plane; Figure (b) show their lifts to 3-space.}
\label{figure:signs of crossings}
\end{figure}

\np There is more to be learned from elementary observations.
Consult Figure \ref{figure:signs of crossings}(a), which shows four little arcs
in the projection of a typical $\cK = \tilde{f}(S^1)$ onto the $xy$ plane.   
The four strands are labeled 1,2,3,4.  First consider strands 1 and 2. Both are
necessarily oriented in the direction of decreasing $x$ because they lie
in the half-space defined by $f'(t) > 0$.  Since $f'$ is decreasing on strand 1,
it follows that $f^{\prime\prime}$ is negative on strand 1, so
$-f^{\prime\prime}$ is positive, so strand 1 lies above the $xy$ plane. Since
$f'$ is increasing on strand 2, it also follows that strand 2 lies below the
$xy$ plane. Thus the crossing associated to the double point in the projection
must be negative, as in the top sketch in Figure \ref{figure:signs of
crossings}(b), and in fact the same will be true for {\em every} crossing in the
upper half-plane. For the same reasons, the projected image of every crossing in
the lower half of the $xy$ plane must come from a positive crossing  in 3-space.
Thus, up to Reidemeister II moves, $\cK$ is a closed braid which factorizes (up
to cyclic permutation) as a product $NP$, where the open braid $N$ (resp. $P$)
represents some number of negative (resp. positive) crossings.  Moreover,  the
type of any such knot is completely defined by its singular projection onto
the $xy$ plane. \\

\np As an example,  the single double point in the left sketch in
Figure \ref{figure:2-braid unknots} lifts to a negative crossing in 3-space,
whereas that in the right sketch lifts to a positive crossing.
Thus $\tilde{f}_+$ defines the 2-braid representative $\sigma_1^{-1}$ of the
unknot and $\tilde{f}_-$ gives the representative $\sigma_1$.\\

\begin{figure}[htpb]
\centerline{\BoxedEPSF{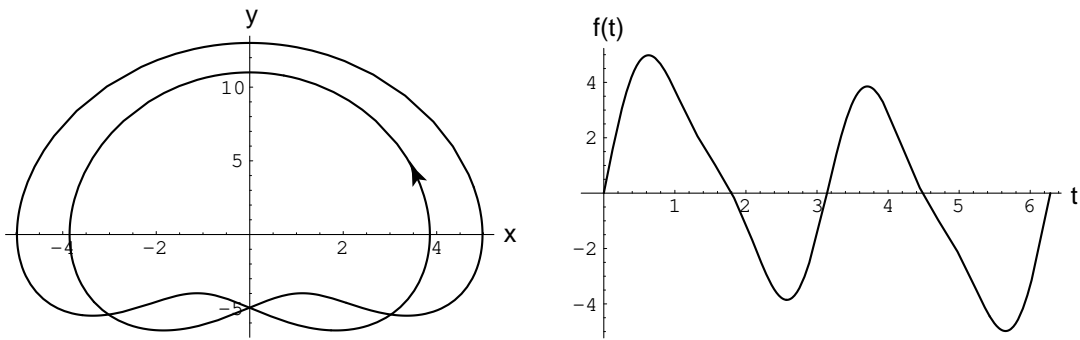 scaled 800}}
\caption{ The function $f(t) = sin(t) + 4 sin(2t) +
sin(4t)$  determines a holonomic positive
$2$-braid trefoil. }  
\label{figure:holonomic trefoil}
\end{figure}

\np A different example is given in Figure \ref{figure:holonomic trefoil}, which 
shows a holonomic parametrization of the 
positive trefoil knot as the 2-braid
$\sigma_1^3 $. The graph of the defining function $f(t)$ is also illustrated. It
has 2 local maxima and 2 local minima and 4 zeros, so the knot is defined by a
holonomic 2-braid. 
The parametrization in Figure \ref{figure:holonomic trefoil} was found
by John Bueti, Michael Kinnally, and Felix Tubiana
during an undergraduate summer research project \footnote {supported by NSF Grant
DMS-98-10750} at Columbia University.  Notice that the graph of $f(t)$
suggests a sawtooth. They were able to generalize this example to 2-braid
representatives of other type $(2,q)$ torus knots by the use of truncated Fourier
approximations of certain sawtooth functions.  However, their partial results
are  a little bit complicated to describe, and at this writing they have not
proved that their functions work for all $q$.  It is clear that much   remains to
be done. \\

\np Note that our holonomic knots are special cases of
Kauffman's {\it Fourier} knots and Trautwein's {\it harmonic} knots, both of which
are  parametrized by three distinct truncated Fourier series in $\reals^3$
\cite{Kauff}, \cite{Traut}.\\ 

\np These knots were introduced into the literature in \cite{Vass}, where they
appeared as a special case of the $n$-jet extension of $f: C^k\to\reals$, where
$C^k$ is the disjoint union of $k$ copies of $S^1$, i.e., the
map $\tilde{f}_n:C^k\to\reals^n$ which is defined by
$\tilde{f}_n(t_1,\dots,t_k) =
(f(t_1,\dots,t_k),f'(t_1,\dots,t_k),\dots,f^{(n-1)}(t_1,\dots,t_k))$.  (Remark:
We have changed Vassiliev's conventions slightly because the 3-jet extension, as
he used it, results in sign conventions which will be confusing for knot
theorists.)  Vassiliev called his $n$-knots {\em holonomic knots} and studied
them.  One of his results for $n=3$ was:

\

\np {\bf Theorem \cite{Vass}:} {\it Every tame knot type in $\reals^3$ can be
represented by a holonomic closed braid of some (in general very high) braid
index.}

\

\np In \cite{Vass} Vassiliev asked the question: ``Is it true that any two
holonomic knots in $\reals^3$ represent the same knot type if and only if they
are isotopic in the space of holonomic knots?"  Figure \ref{figure:non-holonomic
isotopy} gives an example of a holonomic and a non-holonomic isotopy. 
\begin{figure}[htpb]
\centerline{\BoxedEPSF{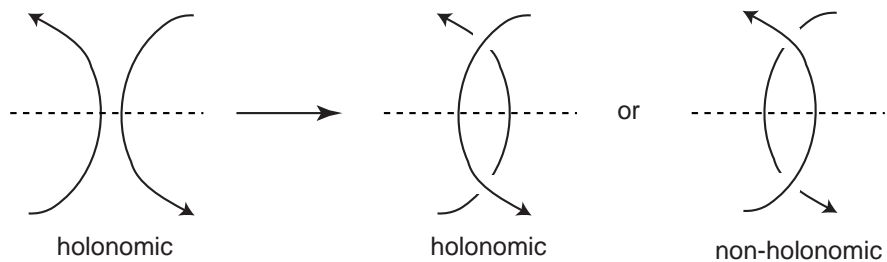 scaled 800}}
\caption{Holonomic and  non-holonomic isotopies on fragments of a knot
diagram.}
\label{figure:non-holonomic isotopy}
\end{figure}
\\

\np The main results in this note are a very simple new proof of a sharpened
version of Vassiliev's Theorem and an affirmative answer to his question: 

\begin{theorem}
\label{theorem:Vassiliev} Every tame knot type in
$\reals^3$ can be represented by a holonomic closed braid $\cH$.  
Moreover the braid index of $\cH$ can be
chosen to be the minimal braid index of the knot type.     
\end{theorem}

\begin{theorem}
\label{theorem:holonomic isotopy}
If two holonomic knots with defining functions $f_0$ and $f_1$ represent the same
knot type, then there is a generic holonomic isotopy $F:S^1\times I\to\reals$
with $F(t,0)=f_0(t)$ and $F(t,1)=f_1(t)$. Thus the study of holonomic knot types
is equivalent to the study of ordinary knot types. 
\end{theorem}

\np Here is an outline of the paper. In $\S$\ref{section:background} we set up
essential background. We will also state and prove Theorem 0, which may be of
interest in its own right. In
$\S$\ref{section:proofs} we prove Theorems 1 and 2.
In $\S$\ref{section:Legendrian} we discuss
parametrizations of Legendrian knots. Since the literature on Legendrian knots
may not be well-known to knot theorists, we will discuss them
in fairly simplistic terms in $\S 4.1$, using insights which we gained as we
struggled to understand the constraints placed by the requirement that
Legendrian knots are tangent to the standard tight contact structure  on
$\reals^3$. Propositions \ref{proposition:holonomic-Legendrian} and 
\ref{proposition:Legendrian cousins} of $\S 4.2$ uncover relationships between
holonomic and Legendrian parametrizations of knots.

\

\np {\bf Remark:} The proofs of the results in this paper make very heavy use of
the work of Garside \cite{Gar} and the related work in {\cite{Adjan},
\cite{Epstein}, \cite{El-M}. Indeed, it appeared to us as we worked out details,
that many of Garside's subtler results seemed to be designed explicitly for
holonomic knots!  We note that this is the first time that we have encountered 
a direct natural connection between the algebraic solution to the word and
conjugacy problems in the braid group, via the work of Garside and others, and
the geometry of knots.  This will become clear after the statement and proof of
Theorem 0.
\\

\np {\bf Acknowledgments:}  We thank Oliver Dasbach, Sergei Chmutov and Serge
Tabachnikov for stimulating discussions and helpful remarks.

\section{Background and notation:} \label{section:background} We
summarize below the main facts we will need from the published literature. The
reader is referred to
\cite{Bir} for background material on the Artin braid groups $\{\bB_n; n
=1,2,\dots \}$. We will use the standard elementary braid generators
$\sigma_1,\dots,\sigma_{n-1}$, where $\sigma_i$ denotes a positive crossing of
the $i^{th}$ and $(i+1)^{st}$ braid strands. Defining relations in $\bB_n$ are: 

\begin{enumerate}
\item [(A1).] $\sigma_i\sigma_j = \sigma_j\sigma_i$ if $|i-j|\geq 2$ for all
$1\leq i,j\leq n-1.$
\item [(A2).] $\sigma_i\sigma_{i+1}\sigma_i = \sigma_{i+1}\sigma_i\sigma_{i+1}$,
where $1\leq i \leq n-1.$
\end{enumerate}

\np We shall use the symbols:

\begin{enumerate}
\item [] $X,Y,Z,\dots $ for words in the generators of $\bB_n$.
\item [] $\bX,\bY,\bZ,\dots$ for the elements they represent in $\bB_n$.
\item [] $\cX,\cY,\cZ,\dots$ for the associated cyclic words.
\item [] $[\cX],[\cY],[\cZ],\dots$ for their conjugacy classes in $\bB_n$.
\item [] $P,P_1,P_2,\dots$ for positive words in the generators of $\bB_n$.
\item [] $N,N_1,N_2,\dots$ for negative words in the generators of $\bB_n$. 
\item [] $X = Y$ if $X$ and $Y$ represent the same element of $\bB_n$.
\item [] $H = NP$ if $H$ can be so represented, with $P$ positive and $N$
negative.
\item [] $\cH = N|P = P|N$ if the cyclic word $\cH$ has such a representation.
\end{enumerate}
\ 

The vertical bar in the notation $\cH = N|P
= P|N$ can be interpreted geometrically as
separating the closed braid into $N$egative (above
the $xy$ plane) and $P$ositive (below the $xy$
plane) pieces. \\

\np {\bf The contributions of Vassiliev.}  We shall use the following
results from \cite{Vass}:

\begin{enumerate}

\item [(V1).] If a knot is represented by a diagram in the $xy$ plane which has
only negative (resp. positive) crossings in the upper (resp. lower) half-plane,
then it may be modified by isotopy to a knot which has a holonomic
parametrization. 

\item [(V2).] Every holonomic knot may be modified by a holonomic isotopy to a
holonomic closed braid, i.e., a closed braid which splits as a product $N|P$.
(Vassiliev calls them {\em normal} braids, but we prefer  the term {\it
holonomic} closed braid.) 

\item [(V3).]  The following modifications in a holonomic closed braid $N|P$
are realized by holonomic isotopy:

\begin{enumerate}
\item Positive (resp. negative) braid equivalences in $P$ (resp. $N$).
\item If $P=N_1P_1$ in $\bB_n$, where $N_1$ is negative and $P_1$ is
positive, replace $N|P$ by $NN_1|P_1$, with a similar move at the other
interface. A special case of this move occurs when we add or delete
$\sigma_j^{-1}\sigma_j$ at the interface. 
\item Insert $\sigma_n^{\pm 1}$ at either interface of the $n$-braid $N|P$ to
obtain an $(n+1)$-braid, or the inverse of this move.    
\end{enumerate}
\end{enumerate}

\np {\bf The contributions of Garside.}  Let $\bB^+_n$ be the semigroup
of positive words in $\bB_n$ which is generated
by $\sigma_1,\dots,\sigma_{n-1}$, with defining relations (A1) and (A2). If
$X,Y$ are words in $\bB^+_n$ we write $X \doteq Y$ to indicate that $X$ and
$Y$ are equivalent words in $\bB^+_n$. In  \cite{Gar} F. Garside proved that the
natural map from $\bB^+_n\to\bB_n$ is an embedding, i.e., if $P_1,P_2$ are
positive words in the generators of $\bB_n$ then $\bP_1 = \bP_2$ if and only if
$P_1 \doteq  P_2$.  The same is true for negative words and negative
equivalences.
\\ 

\np In \cite{Gar} Garside introduced the  $n$-braid ${\bf \Delta}$,
a `half-twist'
which is defined by the word: $$ \Delta = \Delta_n = (\sigma_1\sigma_2\dots
\sigma_{n-1})(\sigma_1\sigma_2\dots\sigma_{n-2})\cdots(\sigma_1\sigma_2)(\sigma_1)$$
and uncovered some of its remarkable properties.  A {\em fragment of} $\Delta$
is any initial subword of one of the (many)  positive words which are
representatives of {\bf $\Delta$}. 

\begin{enumerate}
\item [(G1).] For each $i = 1,2,\dots,{n-1}$ the weak commutativity
relation $\Delta_n \sigma_i \doteq \sigma_{n-i}\Delta_n$ holds.

\item [(G2).] For each $i = 1,2,\dots,{n-1}$ there are fragments of $\Delta$,
say $U_i$ and $V_i$, such that $\Delta \doteq U_i\sigma_i \doteq \sigma_i V_i$,
or equivalently $\sigma_i^{-1} = \Delta^{-1} U_i = V_i\Delta^{-1}$. 

\item [(G3).] (\cite{Gar}, \cite{Adjan}, \cite{Epstein})  For any $\bX\in\bB_n$
and any $X$ which represents $\bX$ there exists a systematic procedure
for converting $X$ to a unique {\em normal form} $\Delta^kP_1\cdots P_r$. In the
normal form each $P_i$ is a fragment of $\Delta$, also each $P_i$ is a longest
possible fragment of $\Delta$ in the class of all positive words which are
positively equal to $P_i$. Finally, $k$ is maximal  and $r$ is simultaneously
minimal for all such representations.  If one starts with 
$X=\Delta^iQ$, where $Q$ is positive, one finds the normal form by repeatedly
`combing out' powers of $\Delta$ from $Q$, i.e. using the fact that if
$k>i$, then $X=\Delta^{i+1}Q_1$ where $Q\doteq\Delta Q_1$. A finite number
of such combings yields $P$. An examination of all positive words which represent
the same element of $\bB_n$ as $P$ produces the decomposition $P\doteq P_1\cdots
P_r$.

\item [(G4).] (\cite{Gar} and \cite{El-M}) For any conjugacy class
$[\cX]\in\bB_n$ and any braid $X$ whose closed braid represents $[\cX]$, there
exists a systematic procedure for converting $X$ to a related normal form
$\Delta^{k'}P'_1\cdots P'_{r'}$ where $k'$ is maximal and $r'$ is minimal for
all such representations of words in the same conjugacy class. Call this a {\em
summit form} for $[\cX]$. The integers $k'$ and $r'$ are unique but the positive
braid $P'_1\cdots P'_{r'}$ is not unique. However there is a finite collection of
all such positive braids and it is unique. 

\item[(G5).] There is a systematic procedure for finding a summit form: Assume
that $X\in\bB_n$ is in the normal form of (G3). Garside proves that there
exists a positive word $W$ which is a product $A_1A_2\cdots A_z$, where each
$A_i$ is a fragment of $\Delta$, such that $W^{-1}XW = X'$, where $X'$
is a summit form.  He also shows how to find $A_1,\dots,A_z$. Let
$W_i=A_1A_2\cdots A_i$. Then each $W_i^{-1}XW_i = 
\Delta^{k_i}P_{1,i}P_{2,i}\cdots P_{r_{i},i}$ where $k_i\geq k_{i-1}$ and 
$r_i\leq r_{i-1}$, also $r=r_1$ and  $r'=r_z$.

\item[(G6).]  If $X' = \Delta^{k'}P'_1\cdots P'_{r'}$ and $X'' =
\Delta^{k'}P''_1\cdots P''_{r'}$ are both summit forms of $X$,
then $X'$ and $X''$ are related by a series of positive conjugacies, as in
(G5), with each $k_i = k'$ and each $r_i = r'$. 
\end{enumerate}

\np We pause in our review of the background material to point out a connection
between Garside's work and holonomic isotopy. 

\

\np {\bf Theorem 0:} {\it Given any open braid $H$ which is in holonomic form $NP$,
the following hold:}
\begin{enumerate}
\item {\it The open braid $NP$ may be brought to Garside's form $\Delta^{-q}Q$,
where $q$ is a non-negative integer and
$Q\in\bB^+_n$, by a holonomic isotopy.}  
\item {\it A further holonomic isotopy converts the open braid $\Delta^{-q}Q$
to its unique normal form $\Delta^kP_1P_2\cdots P_r$.}  
\item {\it A final holonomic isotopy
converts the associated closed braid to one of its summit forms
$\Delta^{k'}P'_1\cdots P'_{r'}$.} 
\end{enumerate}

\np{\bf Proof of Theorem 0:} \begin{enumerate}
\item We are given $H = NP$.  If $N = \emptyset$ there is
nothing to prove, so assume that
$N
\not=
\emptyset$. Using only negative equivalences,  comb out powers of $\Delta^{-1}$
to write $N$ in the form
$N=\Delta^{-k}\sigma_{j_1}^{-1}\cdots\sigma_{j_{s-1}}^{-1}\sigma_{j_s}^{-1}$
where the word $\sigma_{j_1}^{-1}\cdots\sigma_{j_{s-1}}^{-1}\sigma_{j_s}^{-1}$ is
not equivalent to any word which contains a power of $\Delta^{-1}$. By (V3) this
move can be realized by a holonomic isotopy. If
$s=0$ the theorem is true, so assume that $s\geq 1$.  By (G2) we may find a fragment
$U_{j_s}$ of
$\Delta$ such that $\sigma_{j_s}^{-1} =
\Delta^{-1}U_{j_s}$. Therefore by (V3) our closed braid is holonomically
equivalent to
$\Delta^{-k}\sigma_{j_1}^{-1}\cdots\sigma_{j_{s-1}}^{-1}\Delta^{-1}U_{j_s}P$,
which by (G1) is holonomically equivalent to
$\Delta^{-k-1}\sigma_{n-j_1}^{-1}\cdots\sigma_{n-j_{s-1}}^{-1}U_{j_s}P.$
Induction on $s$ completes the proof. 
\item To change $\Delta^{-k}Q$ to normal form for its word
class, first comb out powers of
$\Delta$ from the positive part
$Q$. Then put the new positive part into Garside's
normal form. Both of these steps are achieved by positive
braid equivalences, so by (V3) this part of the work is realizable by a
holonomic isotopy.  
\item Following (G5)
we move the normal form just achieved to a summit form and check that this
move is holonomic.  When we consider conjugates of $H$ in (G5), Garside tells us
that it is enough to work with conjugates $W^{-1}HW$, where $W$ is positive. If
$H$ was holonomic before conjugating by $W$, then the positivity of $W$ keeps
the new braid holonomic after conjugation. \endpf
\end{enumerate}

\np {\bf The contributions of Markov.} \  We will need to use Markov's
Theorem  (see
\cite{Bir} or \cite{Morton}) in our work.  We state it in the form in which it
will be most useful in this paper: \\

\np {\bf Markov's Theorem}  Let $X\in \bB_p$ and $X'\in\bB_q$ be two braids whose
associated closed braids $\cX, \cX'$ define the same oriented knot type. Then there
is a sequence of conjugacy classes in the braid groups $\{\bB_n, \ 1\leq n
<\infty\}$:

\begin{center}
$[\cX] = [\cX_0] \longrightarrow [\cX_1] \longrightarrow [\cX_2]
\longrightarrow \cdots \longrightarrow [\cX_{r-1}] \longrightarrow [\cX_r] =
[\cX']$
\end{center}

\np where $[\cX_j] \subset \bB_{n_j}$, and there are open braid representatives
$X_{j,1}$ for each $1\leq j\leq r$ and $X_{j,2}$  for each $0\leq j\leq r-1$ of
the conjugacy class $[\cX_j]$  such that either:

\begin{enumerate}
\item [(M1).] $n_{j+1} = n_j + 1$ and $X_{j+1,1} = X_{j,2}\sigma_{n_j}^{\pm 1}$
(adding a trivial loop), or
\item [(M2).] $ n_{j+1} = n_j - 1$ and
$X_{j+1,1}\sigma_{n_{j+1}}^{\pm1} = X_{j,2}$  (deleting a trivial loop).\\
\end{enumerate}

\section{The proofs:} \label{section:proofs}
\np {\bf Proof of Theorem 1.}  By a well-known theorem of Alexander (see
\cite{Bir} for example), any knot type may be represented as a closed braid. In
the collection of all closed braid representatives of a given knot type, let us
suppose that $\cK$ is the closure of the braid
$\sigma_{\mu_1}^{\epsilon_1}\cdots\sigma_{\mu_r}^{\epsilon_r}$, where each
$\epsilon_j = \pm 1$.  Use (G2) to replace each $\sigma_{\mu_i}^{-1}$ by
$\Delta^{-1}V_{\mu_i}$. Then use (G1) to push the $\Delta^{-1}$'s to the left. 
This replaces the given representative by a new closed braid $\cH$ which 
is in the desired form $N|P$. But then, by (V1), $\cH$ is isotopic to a holonomic
closed braid. 

\  

\np If we had chosen the braid representative of our knot $\cK$ to have minimum braid
index (which we can do without loss of generality) then the holonomic braid $\cH$
will also have minimum braid index because the changes which we introduced to achieve
$N|P$ form do not change the number of braid strands.  
\endpf

\np {\bf Remark:} This proof differs from Vassiliev's proof in the following
way. We both begin with an arbitrary knot diagram. He modifies the given diagram
by a move which eliminates positive (resp. negative) crossings in the upper
(resp. lower) half-plane, at the expense of adding some number of anticlockwise
loops which encircle points on the $x$-axis. He then changes the resulting
holonomic knot to a holonomic closed braid by using holonomic Reidemeister II
moves, as in our Figure \ref{figure:non-holonomic isotopy}.  In our proof,
we modify the original diagram to a closed braid, and then use a very small part
of Garside's work, without changing the number of braid strands, to complete the
proof. We will see this theme expanded in the proof of Theorem \ref{theorem:holonomic
isotopy}.\\

\np {\bf Proof of Theorem 2.}  We begin our proof of Theorem 2 with two
holonomic knots which, by (V1) and (V2), can be assumed to be holonomic closed
braids. Thus $\cH = N|P$ and $\cH' = N'|P'$. By hypothesis, our closed braids
define the same oriented link types in $\reals^3$.  Markov's Theorem then gives
us a chain of conjugacy classes of braids which connects them:

$$[\cH] = [\cX_0] \longrightarrow [\cX_1] \longrightarrow [\cX_2] \longrightarrow
\cdots \longrightarrow [\cX_{r-1}] \longrightarrow [\cX_r] = [\cH'].$$

\np Let's consider the passage from the class $[\cX_j]$ to the class
$[\cX_{j+1}]$. By Markov's Theorem, we must choose representative open braids
$X_{j,1}, X_{j,2}$ of the conjugacy class $[\cX_j]\subset \bB_{n_j}$ and show
that in either of the two cases (M1), (M2) our representatives and the Markov
moves between them are holonomic.  

\begin{itemize}

\item  In the situation of (M1) the braid $X_{j,2}\in \bB_{n_j}$ is not
necessarily holonomic.  We change the closure of $X_{j,2}$ to a holonomic closed
braid $H_{j,2}$, if necessary, using Garside's methods:
$H_{j,2} = \Delta_{n_j}^{-p}P$  where $p\geq 0$ and $P\in\bB^+_{n_j}$.  Then
$\Delta_{n_j}^{-p}P \sigma_{n_j}^{\pm 1} = \cH_{j+1,1}$ is holonomic  for both
choices of the exponent of $\sigma_{n_j}$.  By (V3) the
passage $\cH_{j,2}\longrightarrow\cH_{j+1,1}$, which adds a trivial loop at
the interface between the positive and negative parts of $\cH_{j,2}$, can be
realized by a holonomic isotopy.

\item  In the situation of (M2) the braid $X_{j+1,1}\in \bB_{n_j-1}$ is
in general not holonomic, but we may replace it as above with a holonomic
$(n_j-1)$-braid $H_{j+1,1}$ which is in the form $N|P$.   But then $\cH_{j,2}$ is
also holonomic, for both choices of the exponent of  $\sigma_{n_j-1}$, because
the ambiguously signed letter is at the interface between the positive and
negative parts of $\cH_{j,2}$.  The passage $\cH_{j,2}\longrightarrow
\cH_{j+1,1}$, which deletes the trivial loop, can clearly be realized by a
holonomic isotopy.
\end{itemize} 

\np The initial and final braids $\cH_{0,1}$ and $\cH_{r,2}$ have not yet been
defined; we take them to be the given representatives $N|P$ and $N'|P'$ of
$[\cH]$ and $[\cH']$ respectively.  \\

\np Thus we have produced a sequence of conjugacy classes of holonomic braids
which joins $[\cH]$ to $[\cH']$, and in each conjugacy class we have two
holonomic representatives $\cH_{j,1}$ and $\cH_{j,2}$, and for each
$j=1,2,\dots,r-1$ we go from $\cH_{j,2}$ to $\cH_{j+1,1}$ via a holonomic
isotopy.
\\

\np The only point which remains to be proved is that, in each conjugacy class
$[\cX_j] \subset \bB_{n_j}$ in the sequence, there is a holonomic
isotopy between the two chosen holonomic representatives.  The proof will be 
seen to be independent of
$j$, so we simplify the notation, setting $n=n_j$. Assume that we have two
holonomic representatives $\cH, \cH'$ of the same conjugacy class in $\bB_n$. 
Our task is to prove that they are holonomically isotopic. 
\\

\np In view of Theorem 0,  we may assume without loss
of generality that $\cH$ and $\cH'$ are summit forms.  By (G6), two summit forms in
the same conjugacy class are related by conjugations by positive words which are
fragments of $\Delta$ using only positive conjugacy.  Again, positive
conjugation by a positive word sends holonomic braids to holonomic braids.  Thus
we have the desired holonomic isotopy from
$\cH$ to
$\cH'$ and the proof of Theorem 2 is complete.
\endpf

\section{Legendrian knots}
\label{section:Legendrian} 
Theorem 2 was in some ways an unexpected result, because it is well-known
that the analogous theorem fails in the case of Legendrian knots. With the goal
of understanding this situation better we investigate, briefly, 
parametrizations of Legendrian knots.  We begin by introducing Legendrian
knots, via their parametrizations. 

\subsection {Parametrizing Legendrian knots} 
\label{subsection:parametrizing Legendrian knots}  Introduce coordinates $(x, v, z)$ in
$\reals^3$.  A knot is {\em Legendrian} if it is represented by a smooth embedding of $S^1
\rightarrow
\reals^3$ whose image is everywhere tangent to the planes of the standard
contact structure on
$\reals^3$. The {\em standard} contact structure on $\reals^3$ is the
nonintegrable field of planes defined by the differential 1-form $\alpha = z dx -
dv$.  A curve will be everywhere tangent to these 2-planes if it is in the
kernel of
$\alpha$, so if a knot can be parametrized by $(x(t), v(t), z(t))$, where
$x,v,z$ are real-valued periodic functions with period (say) $2\pi$, then
it is tangent to this contact structure if and only if $z(t) = dv/dx =
v'(t)/x'(t)$ for all $t\in [0,2\pi]$. \\

\np A simple example is given by the parametrization 
$(-cos(t), -sin^3(t), -3 sin(t) cos(t))$.  See Figure
\ref{figure:Legendrian unknot}.  
In this example the projection onto the $xv$ plane (the so-called {\em front
projection}) has cusps, but they do not represent points of non-smoothness in the
space curve. Indeed, the front projection of a Legendrian knot always contains
$2m$ cusps for some $m\geq 1$.  Notice that  the $z$ coordinate is
the slope of the tangent to the curve in the front
projection. This  makes these projections similar to our $xy$ projections of
holonomic knots, because the signs of the crossings in the front projection
are completely determined without further data. Here is a simple rule for finding the signs of
the crossings: The front projection has no vertical tangencies, because $z(t)$ would be
undefined at a vertical tangency. Therefore at a double point both
branches intersect a vertical line through the double point transversally. Crossings are
positive (respectively negative) if the two branches at a double point intersect a
vertical line through the double point from opposite (respectively the same)
directions. 
\\

\begin{figure}[htpb]
\centerline{\BoxedEPSF{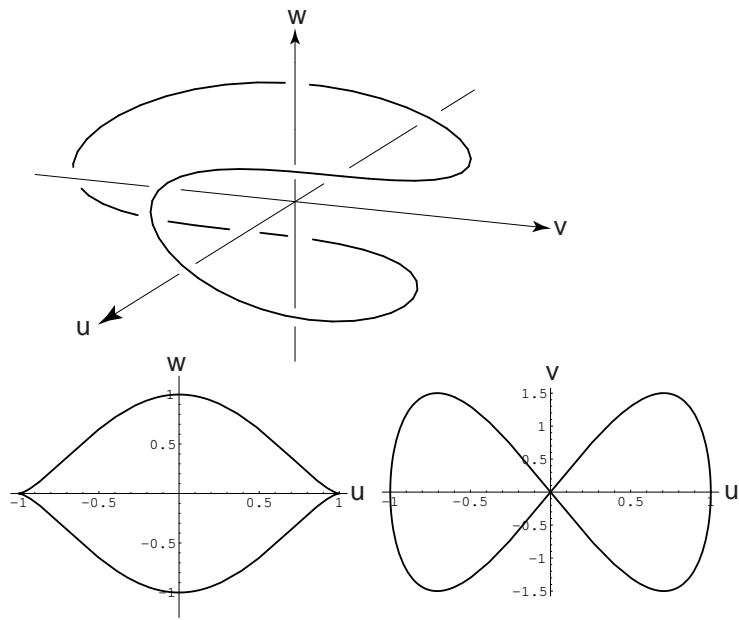 scaled 700}}
\caption {Legendrian unknot $(-cos(t),
-sin^3(t), -3 sin(t) cos(t)).$}
\label{figure:Legendrian unknot}
\end{figure}

\np It is simple to construct other examples (indeed all
possible examples) of parametrized Legendrian knots. Let $x(t)$ and $z(t)$ be
real-valued smooth periodic functions with period $2\pi$.  Then $x(t), z(t)$
have Fourier expansions: 
$$x(t) = \sum_i a_i sin(it) + b_i cos(it), \ \ \ z(t) = \sum_j c_j
sin (jt) + d_j cos(jt).$$
The third coordinate $v(t)$ has a similar expansion, determined up to an
additive constant from the first two by the condition $v'(t) = x'(t) z(t)$. 
Of course we need to be sure that  $x(t)$ and $z(t)$ are generic
(i.e. that they satisfy restrictions analogous to the bulleted assumptions in
$\S$1).

\subsection{The Legendrian cousins of holonomic knots}
\label{section:Legendrian cousins}
Holonomic knots have a parametrization
$(x(t),y(t),z(t))$ where $x(t) = -f(t)$, $y(t) = f'(t)$, and $z(t) =
-f''(t))$, and $f$ is a function chosen as in
section \ref{section:Introduction}.  Following a suggestion of S. Chmutov, we
consider now the differential 1-form to whose kernel this parametrization
corresponds. The holonomic parametrization imposes the following conditions on
the three coordinates: $-y(t) = x'(t)$ and
$-z(t) = y'(t)$.  If we assume
$y(t) \neq 0$, then $z = (y dy)/dx$.  Hence the 1-form $\beta$ whose
kernel contains a holonomic knot is given by $\beta = zdx - ydy$.  For the
1-form $\beta$ to define a contact structure there is an additional
condition it must satisfy: it must have a non-vanishing volume form $
\beta\wedge d\beta$. Here $\beta \wedge d\beta = dx \wedge ydy \wedge dz
$, so $\beta\wedge d\beta$ is only non-zero off the plane $y = 0$, and our
earlier assumption that $y\neq 0$ is an essential part of the story. By
construction our holonomic knots are tangent to this contact structure, and we
show further that this contact structure is isomorphic to the standard one on
$\reals^3$:

\begin{proposition}
\label{proposition:holonomic-Legendrian}
Every holonomic knot is tangent to a contact
structure on $\reals^3 - \{(x,y,z) \in \reals^3 : y = 0\}$ which is isomorphic to the standard
contact structure on $\reals^3$ in the complement of the plane $\{(x, v, z) \in \reals^3 : v
= 0\}$.
\end{proposition}

\pf We map the upper half-space $U_+ = \{(x,y,z) \in \reals^3 : y > 0\}$ to the upper
half-space $V_+ = \{(x,v,z) \in \reals^3 : v = y^2/2>0\}$ and the lower
half-space
$U_-  = \{(x,y,z) \in \reals^3 : y < 0\}$ to $V_+$ with the same transformation $y \rightarrow
y^2/2$.  This gives us transformations sending the two half-spaces of $\reals^3 - \{y = 0\}$
to the upper half-space
$V_+$ of $\reals^3$.  In $V_+$ we have the new relations $dv = ydy$ and $dv = zdx$, so our
1-form $\beta = zdx - ydy$ becomes the 1-form  $\alpha = zdx -
dv$.  The contact structure defined by $\alpha$ on $V_+$ is exactly the restriction of the
standard one on $\reals^3$.  We can extend it from $V_+$ to all of $\reals^3$ then, although
on all of $\reals^3$ the transformation from $\beta$ to $\alpha$ is not invertible.  (It is
not invertible along the plane $v = 0$, where $z = dv/dx$ goes to infinity.)  \endpf

\np In view of Proposition \ref{proposition:holonomic-Legendrian}, a simple
modification of our holonomic parametrizations suggests itself as a method for
parametrizing a related class of Legendrian knots.  Let $\cH$ be any knot type
in $\reals^3$. Using Theorem 1, choose a smooth periodic function $f(t)$ such
that
$(-f(t),f'(t),-f^{\prime\prime}(t))$ is a holonomic parametrization $H(t)$ of
$\cH$. Let
$(-f(t),f'(t)^3,-3f'(t)f^{\prime\prime}(t))$ be a parametrization $L(t)$ of
a new knot $\cL$. We call $\cL$ the {\it
Legendrian cousin} of the holonomic knot $\cH$.  We already gave one example in Figure 7, 
which illustrated the Legendrian cousin of the unknot of Figure 1. A different
example is given in Figure 8, which shows the projection onto the $xy$-plane
of the Legendrian cousin of the holonomic trefoil of Figure 5. 

\begin{figure}[htpb]
\centerline{\BoxedEPSF{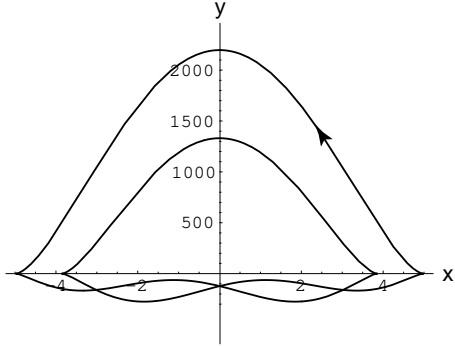 scaled 700}}
\label{figure:cousin of holonomic trefoil}
\caption {Legendrian cousin of the holonomic trefoil of Figure 5.}
\end{figure}
      
\begin{proposition}
\label{proposition:Legendrian cousins}
Let $\cH$ be any knot type in $\reals^3$. As above, let $H(t)$ be its holonomic
parametrization and $L(t)$ the parametrization of $\cL$, its Legendrian
cousin.   Then:
\begin{enumerate}

\item The parametrized curve $L(t)$ is Legendrian, relative to the standard
contact structure $z dx -dv$ on $\reals^3 - \{(x, v, z) \in \reals^3 : v
= 0\}$.

\item The curves $H(t)$ and $L(t)$ are closed
braids with respect to the $z$  axis.  Their
projected images onto the $xy$ plane define
equivalent immersions, i.e., there is a $1-1$
correspondence between the singularities of the
two projections and an isotopy from one projection
to the other induced by the isomorphism of
$\reals^2$ taking coordinates $(x,y) \to (x,v)$.

\item All crossings in the front projection for $\cL$ are negative.

\item Whereas any knot type $\cH$ is represented by a holonomic $H(t)$, the
knot type $\cL$ which is represented by its Legendrian cousin $L(t)$ is always a
fibered knot, i.e. its complement fibers over the circle with fiber a closed
orientable surface $F$, with $\partial F = L$. In fact,
$\cL$ belongs to a very special class of fibered knots which  can be represented
by closed negative braids.
\end{enumerate}
\end{proposition}

\pf To prove (1), notice that
$$-3f'(t)f^{\prime\prime}(t) = {d((f'(t))^3 /dt \over -f'(t)}.$$
To prove (2), let $D(H)$ and $D(L)$ be the parametrized immersed curves
$(-f(t),f'(t))$ \\
and $(-f(t),f'(t)^3)$ in the $xy$ plane so that
$D(H)$ is a holonomic projection for
$\cH$ and $D(L)$ is a front projection for $\cL$.  Then the diagrams $D(L)$ and
$D(H)$ are the same immersions, up to the isomorphism of $\reals^2$ 
which changes coordinates
from $(x, y)$ to $(x, v)$.  For observe that
$$(-f(t_1),f'(t_1)) = (-f(t_2),f'(t_2)) \Longleftrightarrow
(-f(t_1),f'(t_1)^3) = (-f(t_2),f'(t_2)^3).$$
Notice that
the cusps in $D(L)$ occur precisely at the points where $D(H)$ crosses
the
$x$ axis.
The assertion about closed braids is then clear, from our work
in $\S$1 of this paper. \\

\np (3) The statement about the signs of the
crossings in $\cL$ follows from the fact that,
using the rule given in Section
\ref{subsection:parametrizing Legendrian knots},
the signs of the crossings in $D(L)$  are always
negative.  However, in $D(H)$ they are negative in
the upper half-plane and positive in the lower
half-plane.  It's easy to understand why this is
the case. At a double point in both $D(h)$ and
$D(L)$ we have distinct values $t_1,t_2$ with
$f(t_1) = f(t_2)$ and $f'(t_1) = f'(t_2)$. The
signs of the crossings are determined by the
$z$-values on the two branches, that is by
$-f^{\prime\prime}(t_1)$ and
$-f^{\prime\prime}(t_2)$ in $H(t)$, but by
$-f^{\prime\prime}(t_1)f'(t_1)$ and
$-f^{\prime\prime}(t_2)f'(t_2)$ in
$L(t)$. From this it follows that at a double
point in the upper half-plane, where
$f'(t_1)=f'(t_2) > 0$, the signs of the crossings
will be the same in $H(t)$ and $L(t)$, whereas at
a double point in the lower half-plane, where
$f'(t_1)=f'(t_2) < 0$, the signs of the crossings
in $L(t)$ are opposite to those in $H(t)$. From
this it follows that all crossings in
$L(t)$ are negative.  An example can be seen in
Figures 8 and 5.  The holonomic knot in Figure 5 is
a positive type (2,3) torus knot, whereas the
Legendrian knot in Figure 8 is a negative type
(2,3) torus knot.\\

\np (4) Finally, it is well-known (see, for example,
\cite{St}) that knots which can be represented by closed negative or 
positive braids are a
very special subclass of fibered knots.\endpf  

\np Proposition 2 shows that the failure of the contact structure defined
by $\beta$ to extend across the plane
$y=0$ has profound consequences. By Theorem 2 an arbitrary isotopy between
holonomic knots can always be deformed to a holonomic isotopy, but it should now
be clear that we cannot expect any such result for the Legendrian cousins of
holonomic knots. Indeed, our proofs of Theorems 1 and 2 above depended
heavily on our ability to move across the plane $y=0$ via holonomic isotopies,
but for Legendrian cousins we are restricted to Legendrian isotopies in either
half-space, so entirely new methods would be needed to prove an analogous
theorem for Legendrian cousins.  \\

\np {\bf Remark:} The Legendrian cousins of
Proposition
\ref{proposition:Legendrian cousins} can be generalized to an infinite family of
Legendrian cousins of
$H(t)$, with parametrizations 
$$L_k(t) =
(-f(t),f'(t)^{2k+1},-(2k+1)f'(t)^{2k-1}f^{\prime\prime}(t)).$$ 
They share a common knot diagram, so they all represent the same knot type. In
fact, they represent the same Legendrian knot type. An explicit
Legendrian isotopy $L(t,s)$ from
$L_k(t)\to L_m(t)$ was suggested by Oliver Dasbach. It begins at
$L(t,0) = L_k(t)$ and ends at 
$L(t,1) = L_m(t)$ and is defined by $L(t,s) =$
$$
(-f(t),sf'(t)^{2k+1}+(1-s)f'(t)^{2m+1},-s(2k+1)f'(t)^{2k-1}f^{\prime\prime}(t)
-(1-s)(2k+1)f'(t)^{2m-1}f^{\prime\prime}(t))$$

\np \small {Joan S. Birman, Math. Dept, Barnard College of Columbia
University,}
\np \small {604 Mathematics, 2990 Broadway, New York, N.Y. 10027.}\\

\np \small {Nancy C. Wrinkle, Math. Dept., Columbia University,}
\np \small {408 Mathematics, 2990 Broadway, New York, N.Y. 10027.}

\end{document}